\documentclass{amsart}
\usepackage{epsfig}

\begin{document}

\newtheorem{thm}{Theorem}[subsection]
\newtheorem{lem}[thm]{Lemma}
\newtheorem{cor}[thm]{Corollary}
\newtheorem{conj}[thm]{Conjecture}
\newtheorem{qn}[thm]{Question}

\theoremstyle{definition}
\newtheorem{defn}[thm]{Definition}

\theoremstyle{remark}
\newtheorem{rmk}[thm]{Remark}
\newtheorem{exa}[thm]{Example}

\def\square{\hfill${\vcenter{\vbox{\hrule height.4pt \hbox{\vrule width.4pt
height7pt \kern7pt \vrule width.4pt} \hrule height.4pt}}}$}

\def\R{\mathbb R}
\def\Z{\mathbb Z}
\def\CP{\mathbb {CP}}
\def\H{\mathbb H}
\def\E{\mathbb E}
\def\F{\mathcal F}
\def\D{\mathcal D}
\def\C{\mathbb C}
\def\til{\widetilde}
\def\N{\mathbb N}
\def\T{\mathcal T}
\def\G{\mathcal G}
\def\S{\mathcal S}
\def\a{{\text{amb}}}
\def\supp{{\text{supp}}}
\def\homeo{{\text{Homeo}}}

\newenvironment{pf}{{\it Proof:}\quad}{\square \vskip 12pt}

\title{The Gromov norm and foliations}

\author{Danny Calegari}
\address{Department of Mathematics \\ Harvard \\ Cambridge, MA 02138}
\email{dannyc@math.harvard.edu}

\maketitle

\begin{abstract}
We define a norm on the homology of a foliated manifold, which refines
and majorizes the usual Gromov norm on homology. This norm depends
in an upper semi--continuous way on the underlying foliation, in the
geometric topology. We show that this norm is non--trivial --- i.e. it distinguishes
certain taut foliations of a given hyperbolic $3$--manifold.

Using a homotopy--theoretic refinement, we show that a taut foliation whose leaf space 
branches in at most one direction cannot be the
geometric limit of a sequence of isotopies of a fixed taut foliation whose
leaf space branches in both directions.
Our technology also lets us produce examples of taut foliations which cannot
be made transverse to certain geodesic triangulations of hyperbolic
$3$--manifolds, even after passing to a finite cover.

Finally, our norm can be extended to actions of fundamental groups of
manifolds on order trees, where it has similar upper semi--continuity
properties.
\end{abstract}

\section{Introduction}

The study of group actions on trees and tree--like objects
has for a long time been an important tool in $3$--manifold topology.
J. Stallings pioneered this approach with a topological proof of
Grushko's theorem \cite{jS65}, and more generally it is by now a standard
observation that a decomposition of a $3$--manifold along an 
incompressible surface is ``dual'' to some action of $\pi_1(M)$ on 
a tree. M. Culler and P. Shalen \cite{CS83} used the action of
$\pi_1(M)$ for $M$ a hyperbolic manifold on the Bruhat--Tits tree of
$SL(2,F)$, where $F$ is the function field of a curve in the
$SL(2,\C)$ representation variety of $\pi_1$, to obtain striking
topological results about $M$. More generally, a sequence of hyperbolic
(or merely negatively curved) structures on a fixed manifold which are
not precompact in the Gromov--Hausdorff topology may be rescaled and filtered
to give in the limit an action of $\pi_1(M)$ on an $\R$--tree (see
e.g. F. Paulin \cite{fP88}). In a $2$--dimensional setting, this idea
is implicit in Thurston's compactification of Teichm\" uller space by
projective measured laminations \cite{wT80}.

However, the tree--like structures on which $\pi_1(M)$ acts in all these
cases admit some kind of invariant measure structure. For a
taut foliation or an essential lamination, the existence of such a
transverse structure is rare, and leads to strong topological conditions on
the underlying manifold. Consequently many ``naturally occurring'' actions
of fundamental groups of $3$--manifolds on non--Hausdorff simply connected
$1$--manifolds and order trees admit no invariant measure. Nevertheless, one
would like to quantify the amount of branching of such trees in some
natural way. In this paper, we introduce a norm on the homology of a
foliated manifold, which is a refinement of the usual Gromov norm on
homology, where one restricts the admissible chains representing a 
homology class to those which are {\em transverse} --- that is, each
singular map $\sigma:\Delta^i \to M$
in the support of an admissible chain must induce a standard foliation on
$\Delta^i$, one which is topologically conjugate to an affine foliation.
For a hyperbolic manifold $M^n$, the Gromov norm of $[M]$ is proportional to
the volume of $M$, and for $n\ge 3$,
a chain whose norm is close to the infimum actually
``detects'' the geometry of $\til{M} = \H^n$ (this is just a restatement of
Mostow's rigidity theorem). For a foliated manifold, the tension between
the geometry of the ambient manifold and the local affine structure 
determined by the foliation can be used to show that the foliated norm
differs from the usual norm in certain cases, which reflect the topology
and the geometry of the foliation.

In particular, we have the following theorems:

\vskip 12pt

{\noindent
{\bf Theorem 2.2.10.}
{\it Let $\F$ be a foliation of $M^n$ whose universal cover is topologically
conjugate to the standard foliation of $\R^n$ by horizontal $\R^{n-1}$'s. Then
$$\|[M]\|_G = \|([M],\F)\|_{FG}$$}}

\vskip 12pt

{\noindent
{\bf Theorem 2.5.9.}
{\it Suppose that $\F$ is a taut foliation with one--sided branching. Then there
is an equality of norms
$$\|[M]\|_G = \|([M],\F)\|_{FG}$$}}

\vskip 12pt

{\noindent 
{\bf Theorem 2.4.5.}
{\it Suppose $M^n$ is hyperbolic and $\F$ is asymptotically separated. Then
$$\|[M]\|_G < \|([M],\F)\|_{FG}$$}}

Here we say that a foliation $\F$ of a hyperbolic manifold
is {\em asymptotically separated} if for
some leaf $\lambda$ of $\til{\F}$, there are a pair of open hemispheres
$H^+,H^- \subset \H^n$ in the complement of $\lambda$ which are separated by $\lambda$.
We point out that a standard conjecture would imply that a taut foliation of a
hyperbolic $3$--manifold is asymptotically separated iff $\til{\F}$ has two--sided branching.

It should be mentioned that a norm for {\em foliations with transverse 
measures} was defined by Connes
(unpublished) and developed in \cite{mG91} and \cite{mG99}. This norm
uses generalized simplices which are simplicial in the tangential
direction and measure--theoretical in the transverse direction. It has
the usual proportionality properties for foliations whose leaves are all
locally isometric to a space of constant curvature. The ``fundamental cycle''
on which this norm is evaluated is really attached to the measured
foliation, and not to the underlying manifold. 
By contrast, our definition is closer in spirit to norms for stratified 
or decorated spaces.

With our definition, the norm on a homology class is 
upper semi--continuous as a function of
the underlying foliation, in the geometric topology. Since the norm
is defined topologically, this gives obstructions for the existence of
a family of isotopies of a fixed topological foliation to converge 
geometrically to some other foliation. The leaf space of the universal
cover of a taut foliation is a (typically non--Hausdorff) simply--connected
$1$--manifold. The non--Hausdorffness comes from branching of the leaf space.
This branching can occur in both directions, in only a single direction,
or not at all. The taut foliation is said in these three cases {\em to have
branching in both directions}, {\em to have branching in only one direction}, 
or {\em to be $\R$--covered}. We show 

\vskip 12pt
{\noindent
{\bf Corollary 3.1.5.}
{\it Let $\F$ with branching on at most one side and $\G$ with two--sided branching, be
taut foliations of $M^3$. Then there is no sequence of isotopies $\G_i$ of $\G$
which converges geometrically to $\F$.}}
\vskip 12pt

Using similar techniques, 
we show 

\vskip 12pt
{\noindent
{\bf Theorem 3.3.2.}
{\it Let $M$ be a hyperbolic $3$--manifold, and $\F$} any {\it taut foliation
with $2$--sided branching. Then there is a geodesic triangulation $\tau$ of $M$
which cannot be made transverse to $\F$. Furthermore, $\tau$ cannot be made
transverse to $\F$ in any finite cover (i.e. $\tau$ is not virtually fine).}}

\vskip 12pt

Problem 3.16 in Kirby's problem book \cite{rK97} asks for a reasonable
real--valued function on the set of $3$--manifolds which measures the
complexity of $\pi_1(M)$ and behaves appropriately under finite covers
and positive degree maps. One may translate this problem into the
foliated context, where one considers pairs $(M,\F)$, finite covers
and {\em transverse} positive degree maps (a map $f:(M,\F) \to (N,\G)$ 
between foliated manifolds is {\em transverse} if every transversal to 
$\F$ is mapped to a transversal to $\G$), in which context our norm
seems like a ``reasonable'' solution.

\vskip 12pt

I would like to thank I. Agol, A. Casson and W. Thurston with whom I had
some interesting discussions about this material. In particular,
I. Agol's work \cite{iA99} 
on volumes of hyperbolic $3$--manifolds with boundary was
particularly inspiring. Furthermore, I
received partial support from a Sloan Dissertation Fellowship and
from the Clay Mathematical Institute while carrying out work on
this paper.

\section{Foliated norms}

\subsection{Foliations}

We give the basic definitions of various kinds of foliations of
$3$--manifolds. More details are to be found in \cite{dG83}.

\begin{defn}
Let $M$ be the subspace of $\R^3$ for which $z\ge 0$, minus the origin.
$M$ has a foliation $\til{\F}$ by leaves of the form $z=\text{const.}$
which are all planes, except for the leaf $z=0$ which is a punctured plane.
This foliation is preserved by the dilation $(x,y,z) \to (2x,2y,2z)$
and so descends to a foliation of the solid torus. This is called the
{\em Reeb foliation} of the solid torus.
\end{defn}

\begin{defn}
A codimension $1$ foliation of a $3$--manifold is {\em Reebless} if it
has no solid torus subsets foliated with a Reeb foliation.
\end{defn}

For Reebless foliations, every leaf is incompressible and the ambient
manifold is irreducible or covered by $S^2 \times S^1$ foliated by horizontal
spheres. Moreover, every loop transverse to a Reebless foliation is
homotopically essential. It follows that a Reebless foliation of a manifold
pulls back in a finite covering to a {\em co--orientable foliation}, one
for which there is a choice of orientation on transversals which is invariant
under leaf--preserving isotopy. Equivalently, there is a $\pi_1(M)$--invariant
orientation on the leaf space in the universal cover.

\begin{defn}
A codimension $1$ foliation of a $3$--manifold is {\em taut} if there is
a circle in the manifold transverse to the foliation which intersects every
leaf.
\end{defn}

Every taut foliation is Reebless. The induced foliation $\til{\F}$ of
the universal cover $\til{M}$ of a Reebless foliation is (topologically)
a foliation of $\R^3$ by planes, which is the product of a foliation of
$\R^2$ by lines with $\R$. Consequently, every plane is properly embedded
and separates $\R^3$ into two topological half--spaces.

\begin{figure}[h]
\scalebox{.5}{\includegraphics{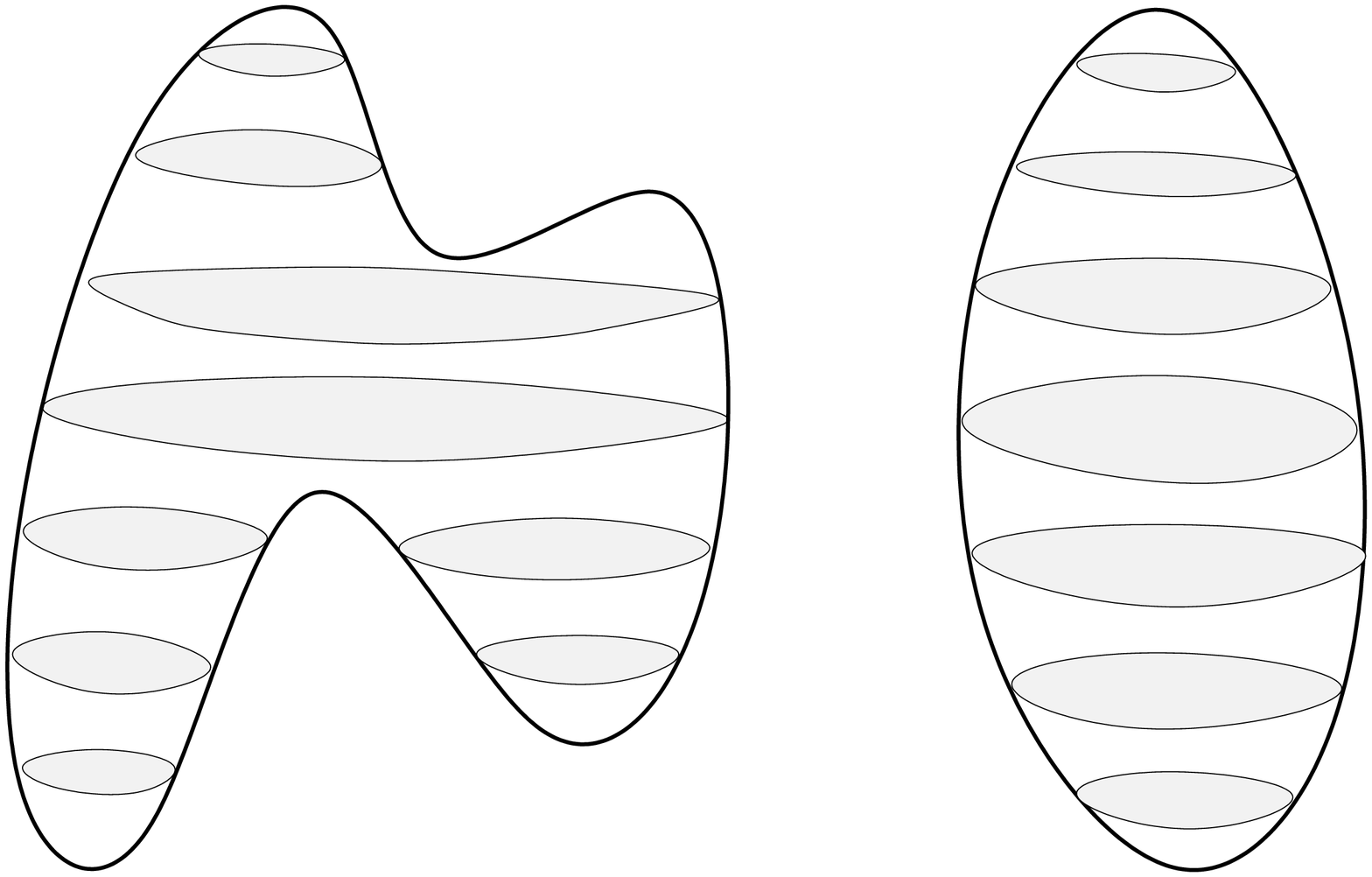}}
\caption{The foliation in the universal cover of two foliations: one
$\R$--covered, the other not.}
\end{figure}

The leaf space $L$ of the universal cover of a taut foliation is a
(possibly non--Hausdorff) simply connected $1$--manifold. The orientation
on this leaf space induces a {\em partial order} on the leaves: 
$\lambda > \mu$ iff there is a positively oriented transversal in
$\til{M}$ from $\mu$ to $\lambda$. The absence of loops in the leaf space
make this partial order well--defined. For readers unfamiliar with the
topology of non--Hausdorff $1$--manifolds, consult \cite{dGwK97}. The
``non--Hausdorffness'' arises from branching of the
leaf--space: an embedded half--open arc in the leaf space 
might have a countably infinite collection of limiting endpoints.
Moreover, this branching might take place at a dense set of points.

\begin{defn}
A taut foliation is $\R$--covered if the pulled back foliation of the
universal cover is topologically conjugate to
the standard foliation of $\R^3$ by horizontal planes.
\end{defn}

For $M$ not finitely covered by $S^2 \times S^1$, the 
leaf space of $\til{\F}$ is Hausdorff exactly when $\F$ is $\R$--covered, for
taut $\F$.

\subsection{Gromov norms}

\begin{defn}
A {\em singular $n$ chain} 
in $M$ is a finite $\R$--linear combination of singular $n$--simplices, 
where a singular $n$--simplex is a map $\sigma:\Delta^n \to M$ from the
standard affine $n$--simplex into $M$. The {\em support} of a chain,
denoted $\supp(C)$ is the set of singular $n$--simplices with non--zero
coefficients in $C$.
\end{defn}

Notice that our convention is to assume that the coefficients in our chains
are allowed to be in $\R$. We assume this without comment in the sequel.

The classical Gromov norm is defined in \cite{mG82} as follows:
\begin{defn}
For $M$ an orientable $n$--manifold, let $[M]$ denote the fundamental
class of $M$ in $H_n(M;\R)$. The {\em Gromov norm} of $M$ is the infimum
of the $L_1$ norm on the singular cycles representing $[M]$. That is,
$$\| [M]\|_G = \inf_{[\sum_i r_i \sigma_i] = [M]} \sum_i |r_i|$$
\end{defn}

\begin{defn}
For $M$ an orientable $n$--manifold and $\F$ a codimension $m$ foliation,
call a singular cycle $\sum_i r_i \sigma_i$ {\em transverse} if the
foliation on the $n$--simplex $\Delta^n$ induced from each singular map
$\sigma_i:\Delta^n \to M$ by pulling back $\F$ is topologically conjugate
to some affine foliation of $\Delta^n$: that is, the foliation by preimages
of points obtained from some affine map $\Delta^n \to \R^{n-m}$. The
{\em foliated Gromov norm} of $[M]$ with respect to $\F$ is defined to be
$$\| ([M],\F)\|_{FG} = \inf_{[\sum_i r_i \sigma_i] = [M]; 
\sigma_i \text{ transverse}} \sum_i |r_i|$$
\end{defn}

\begin{figure}[h]
\scalebox{.5}{\includegraphics{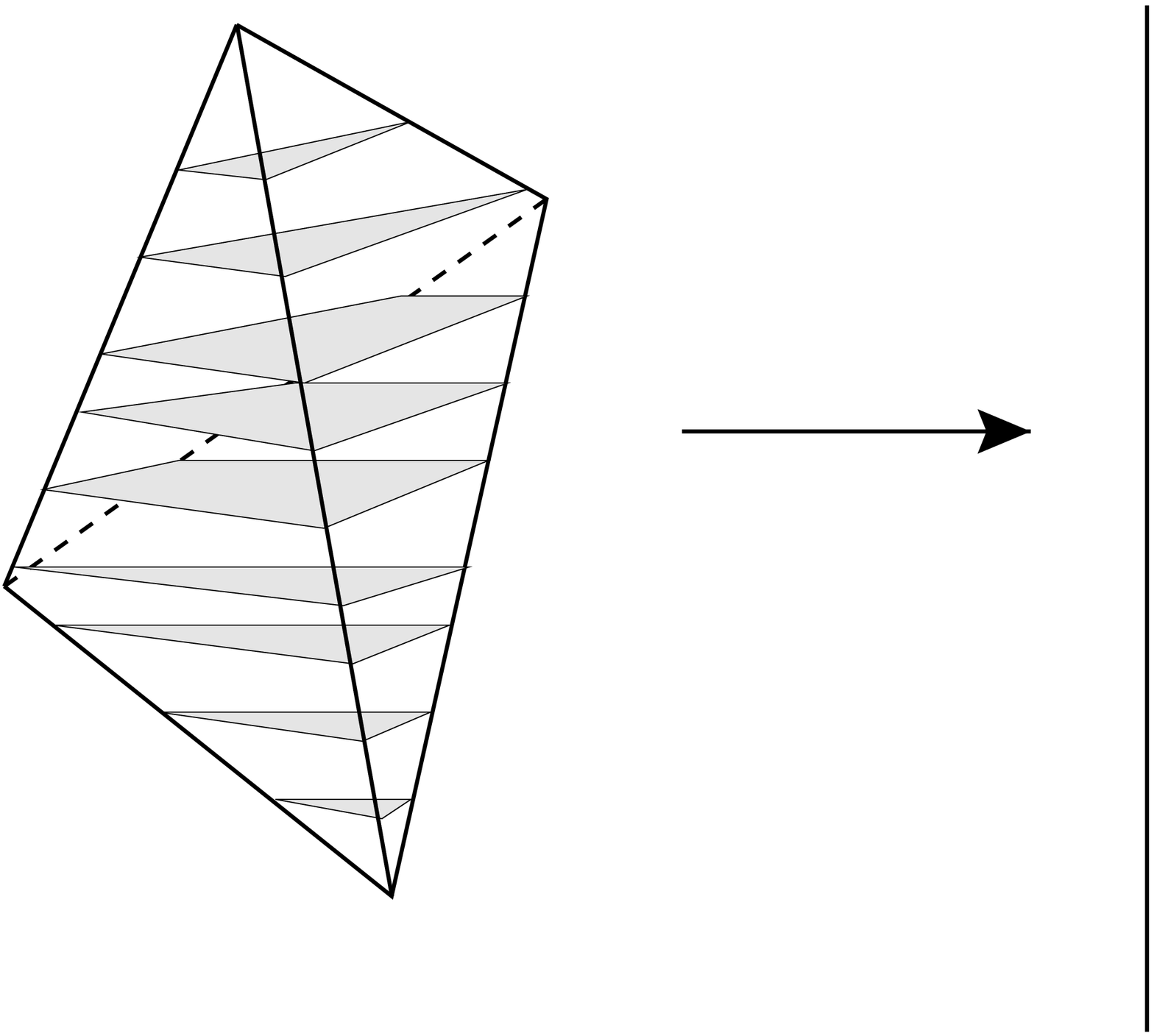}}
\caption{A map from the vertices of a simplex to $\R$ can be extended to
an affine map, which pulls back to an affine foliation of $\Delta^n$.}
\end{figure}

\begin{rmk}
Notice that any map from the vertices of a simplex to $\R$ can be extended
to an affine map of the simplex to $\R$.
\end{rmk}

\begin{thm}
Let $M$ be a $3$--manifold and let $\F$ be a taut foliation of $M$.
\begin{enumerate}
\item{$\| ([M],\F)\|_{FG} < \infty$.}
\item{$\| [M]\|_G \le \| ([M],\F)\|_{FG}$.}
\item{There is a $K(M) < \infty$ such that for any taut foliation $\G$ of
$M$, $\| ([M],\G)\|_{FG} < K(M)$. If $K(M)$ is the infimum of such, define
$\| [M]\|_{FG} = K(M)$.}
\end{enumerate}
\end{thm}
\begin{pf}
Fact 2. is obvious from the definition of the norms. Fact 1. follows from
the existence of a triangulation of $M$ in which the foliation $\F$ can
be put into normal form. Fact 3. follows from the stronger theorem of
D. Gabai that on any $3$--manifold there is a triangulation with respect 
to which {\em every} taut foliation of $M$ can be 
put in normal form. (\cite{dG99})
\end{pf}

\begin{rmk}
Of course, one can define any number of norms on homology by restricting the
class of singular maps which are deemed admissible. The particular restriction
of transversality seems suitable for studying foliations, since it behaves
nicely with respect to many of the usual constructions of foliations, 
e.g. branched covers.
\end{rmk}

\begin{lem}
Let $f:N^n \to M^n$ be a degree $d$ 
branched cover where the branch locus $\gamma$ is
transverse to a foliation $\F$ of $M^n$. Let $\G$ be the foliation obtained
from $\F$ by pullback. Then
$$\| ([N],\G)\|_{FG} \ge d\| ([M],\F)\|_{FG}$$
\end{lem}
\begin{pf}
Let $C = \sum_i r_i \sigma_i$ be a transverse chain representing $[N]$. Then
$f_*C = \sum_i r_i (f \circ \sigma_i)$ is a transverse chain representing
$d[M]$.
\end{pf}

It is a well--known fact that the simplex is distinguished amongst all
affine polyhedra by the fact that any total ordering of its vertices is 
induced from an affine map of the entire simplex to $\R$. In the context
of foliations, this fact has the following generalization:
 
\begin{lem}\label{endpoints_ordered}
Let $\F$ be a foliation of $M^n$ and let $C = \sum_i r_i \sigma_i$ be
a cycle representing $[M]$. Suppose that the leaf space of $\til{\F}$ is
an acyclic $1$--manifold and that the leaves of
$\til{\F}$ are all $\R^{n-1}$. Suppose further for each $i$ that 
$\sigma_i$ lifts to $\til{\sigma_i}:\Delta^n \to \til{M}$ such that the
images of the vertices of $\Delta^n$
inherit a total order from the natural partial order on $L$. Then
$C$ can be ``straightened'' to $C^s = \sum_i r_i \sigma_i^s$ where
$\til{\sigma_i}$ and $\til{\sigma_i^s}$ have the same endpoints in $\til{M}$,
and $C^s$ is a transverse cycle.
\end{lem}
\begin{pf}
Let $S_n = \supp(C)$ and inductively define $S_i = \supp(\partial S_{i-1})$.
Lift each $\sigma \in S_i$ to $\til{\sigma}:\Delta^i \to \til{M}$ and call the
union of some choice of lifts $\til{S_i}$.

Each $\til{\sigma} \in \til{S_1}$ maps its vertices to leaves in $\til{\F}$
which inherit a total order from the partial order on $L$. (Here we actually
require the images of distinct vertices of $\Delta^n$ to lie on distinct
leaves of $\til{\F}$).
It follows we can replace
each $\sigma \in S_1$ by a transverse $\sigma^s$ with the same endpoints.

Let $\partial\sigma:\partial\Delta^i \to \til{\F}$ be transverse. Then
the induced foliation of $\partial\Delta^i$ is the standard foliation of
the unit sphere $S^{i-1}$ in $\R^i$ by its intersection with horizontal
planes. The leaves of this foliation are $(i-2)$--spheres away from the top
and bottom vertex, with respect to the partial ordering on $L$. 
Since leaves of $\til{\F}$ are just $\R^{n-1}$'s, this family of maps of
$(i-2)$--spheres extends to a family of maps of $(i-1)$--balls converging at
the top and bottom to the image of the top and bottom vertices respectively.
This family of maps gives a transverse map $\sigma:\Delta^i \to \til{\F}$ 
agreeing with $\partial\sigma$ on $\partial\Delta_i$.

Thus the straightening procedure can be performed inductively, as required. 
\end{pf}

\begin{rmk}
The ``nondegeneracy'' assumption --- that distinct vertices of a simplex
get mapped to distinct leaves of $\til{\F}$ under lifts of
singular maps $\til{\sigma}$ in the support of $C$ is not really necessary,
since our definition of ``transverse'' cycle includes degenerate affine maps.
In any case, any finite chain can be perturbed chain homotopically to a nondegenerate
chain without violating the total ordering assumption.
\end{rmk}

A taut foliation of $M^3$ has the properties required by 
lemma~\ref{endpoints_ordered}. 

\begin{thm}\label{R_covered}
Let $\F$ be a foliation of $M^n$ whose universal cover is
topologically conjugate to the standard foliation of $\R^n$ by horizontal
$\R^{n-1}$'s. Then $$\| [M]\|_G = \| ([M],\F)\|_{FG}$$
\end{thm}
\begin{pf}
Since $\til{\F}$ is the standard foliation of $\R^n$ by horizontal
planes, the leaf space of $\til{\F}$ is totally ordered. As remarked earlier,
we can perturb a chain $C$ by a chain homotopy 
to a nearby chain $C'$ so that for each $\sigma$
in the support of $C$, $\til{\sigma}$ maps the vertices of $\Delta^n$ to distinct
leaves of $\til{\F}$. It follows from
lemma~\ref{endpoints_ordered} that any chain can be straightened to a
transverse chain.
\end{pf}

In particular, this result holds for $\F$ an $\R$--covered taut foliation
of some $3$--manifold $M$.

\begin{rmk}
One easily extends this argument to see that the foliated Gromov norm
agrees with the usual Gromov norm for foliations whose universal covering
foliations are standard foliations of $\R^n$ by horizontal $\R^{n-m}$'s
--- that is, for {\em product--covered foliations}.
\end{rmk}

\subsection{Measurable chains and equivariant straightening}

Theorem~\ref{R_covered} says that finite chains can be
straightened with respect to an $\R$--covered foliation. An interesting
question is whether the same is true of {\em infinite chains}. We make
this question more precise.

\begin{defn}
Let $\sigma:\Delta^i \to M$ lift to $\til{\sigma}:\Delta^i \to \til{M}$
which can be projected to $\tau:\Delta^i \to L$, the leafspace of $\til{\F}$.
$\sigma$ is {\em monotone} if the stratification of $\Delta^i$ by preimages
of points in $L$ is homotopy equivalent to an affine foliation of $\Delta^i$.
If $\sigma$ is monotone, every $\sigma'$ in the support of $\partial \sigma$
is also monotone.
\end{defn}

\begin{defn}
For each $i$, let $\Sigma_i$ denote the space of singular maps
$\sigma:\Delta^i \to M$ with the compact--open topology. Let
$\Sigma^t_i$ denote the subspace of transverse singular maps.
Let $\Sigma^m_i$ denote the subspace of monotone singular maps.
\end{defn}

In \cite{dC99} we prove the following result:
\begin{thm}\label{product}
Let $\F$ be an $\R$--covered foliation of an atoroidal $3$--manifold $M$.
Then there are co--ordinates on the universal cover $M = \R^2 \times \R$
such that leaves of $\til{\F}$ are horizontal planes $\R^2 \times \text{const}$
and $\pi_1(M)$ acts by elements of $\homeo(\R^2) \times \homeo(\R)$.
\end{thm}

Using this structure,
one has at least the following partial result:

\begin{thm}
Let $\F$ be a co--oriented 
$\R$--covered foliation of a negatively curved closed $3$--manifold $M$. 
Then there are continuous projections from $s_i:\Sigma_i \to \Sigma^m_i$, 
which are compatible with $\partial$ in the sense that
$s_{i-1} \partial = \partial s_i$.
\end{thm}
\begin{pf}
To avoid cumbersome notation, we define instead straightenings of
simplicial maps to $\til{M}$ which are continuous in the compact--open
topology, and which are equivariant under the action of 
$\pi_1(M)$.

For each $\sigma:\Delta^i \to \til{M}$ we can define
two functions $\rho,\tau$ 
in terms of the co--ordinates on $\til{M} = \R^2 \times \R$ by writing
$$\sigma:t \to (\rho(t),\tau(t)) \in \R^2 \times \R$$

\vskip 12pt

We begin by defining $s_1$.
Let $\tau(0) = l, \tau(1) = r$ and let $J$ be the set of numbers bounded
by $r$ and $l$. There is an obvious retract 
$\phi:\R \to J$ which sends everything
above $J$ to the maximum of $J$, and everything below $J$ to the minimum of
$J$. Set $\tau_1(t) = \phi \tau(t)$. Now set
$$\tau'(t) = \inf(s:\text{there exists } t_1\le t \le t_2 \text{ with }
\tau_1(t_1) = \tau_1(t_2) = s)$$
Then $\tau':I \to \R$ is a monotone map, and we can replace the
map $\sigma:t \to (\rho(t),\tau(t))$ with 
$s_1(\sigma):t \to (\rho(t),\tau'(t))$.

\vskip 12pt

Now for $\sigma:\Delta^2 \to \til{M}$, 
we first straighten $\partial \sigma$ using
$s_1$.
The maps $\tau:\partial \Delta^2 \to \R$ are already monotone, so for each
$p \in \Delta^2$, there is a unique $t \in \R$ such that $p$ is in the
convex hull of the points in $\partial \Delta$ which are mapped by $\tau$
to $t$. Define $\tau'(p)$ to be this value $t$, and set
$s_2(\sigma):p \to (\rho(p),\tau'(p))$.

\vskip 12pt

Finally, for $\sigma:\Delta^3 \to \til{M}$, straighten $\partial \sigma$
using $s_2$.
For each $t$ in the interior of the image of
$\tau|_{\partial \Delta^3}$, 
the level set $\tau^{-1}(t) \cap \partial \Delta^3$ is a cellular subset
homotopy equivalent to a circle. It has two frontiers in $\partial \Delta^3$,
on the side where $\tau$ is greater than $t$ and the side where $\tau$ is
less than $t$. Define $D_t \subset \Delta^3$ to be the minimal surface
spanned by the upper frontier of $\tau^{-1}(t)$ (see e.g. \cite{jH86} for
basic facts about minimal surfaces in $3$--manifolds). For distinct $s,t$ the
disks $D_t,D_s$ are disjoint, so we can define ${\tau'}^{-1}(t)$ 
on $\Delta^3$ to be the subset of points above $D_s$ for all $s<t$ and
contained in or below $D_t$. Then set $s_3(\sigma):p \to (\rho(p),\tau'(p))$.

These constructions use only the order structure of $\R$, 
and are therefore equivariant under the action of
$\pi_1(M)$.
\end{pf}

\subsection{The norm is non--trivial on $[M]$}

\vskip 12pt

For $M^n$ hyperbolic, we know $\| [M]\|_G = \text{vol}(M)/v_n$ where $v_n$ is
the volume of the regular ideal $n$--simplex. For a hyperbolic manifold,
any cycle can be chain homotoped to a geodesic cycle by replacing each
singular map of a simplex $\sigma_i:\Delta^n \to M$ with the geodesic simplex
$\sigma_i^g:\Delta^n \to M$ having the same endpoints.

Let $C_j$ be a sequence of geodesic chains 
whose norms converge to the Gromov norm of $M$. Let $\til{C_j}$ denote the
$\pi_1(M)$--equivariant infinite chain obtained by lifting $C_j$ to $\til{M}$.
Let $X$ be the infinite $(n+1)$--valent tree with basepoint, and let
$T$ be the abstract complex obtained by gluing together infinitely many
ideal $n$--simplices along their faces in the pattern described by $X$.
Fix some regular ideal simplex $\Delta \in \H^n$. Then choosing an identification
of $\Delta$ with some simplex of $T$, there is a natural developing map
$\text{dev}:T \to \H^n$ taking each simplex of $T$ to a regular ideal simplex.
If $n=3$ this has as its image the standard regular tessellation of $\H^3$ by
ideal simplices. Otherwise, the representation $Aut(T) \to Isom^+(\H^n)$ is
indiscrete.

\begin{lem}
With notation as above, for any $t, \epsilon > 0$ there is a $j$ such that
for any $k>j$, there is a collection $S_k$ of singular maps in the support
of $\til{C_k}$ and an element $\alpha_k \in Isom^+(\H^n)$ such that $S_k$ and
$\alpha_k\text{dev}(T)$ agree on the ball of radius $t$ about $0$ to within $\epsilon$.
\end{lem}
\begin{pf}
Let $C_k = \sum_i r_i \sigma_i$ be a geodesic chain which very nearly realizes the
Gromov norm of $M$. Then by definition $$\sum_i |r_i| < \text{vol}(M)/v_n + \delta$$
for some small $\delta$. On the other hand,
$$\sum_i r_i \text{vol}(\sigma_i(\Delta^n)) = \text{vol}(M)$$
so the weighted average
$$\frac {\sum_i r_i \text{vol}(\sigma_i(\Delta^n))} {\sum_i |r_i|} > v_n - \delta'$$
for some small $\delta'$; that is, ``most'' of the $\sigma_i(\Delta^n)$, as weighted
by $r_i$, have volume very close to
$v_n$. This implies that they are geometrically very close to regular ideal
simplices, on a big compact set containing most of their mass. (see e.g.
\cite{rBcP92})

Fix a fundamental domain $D$ in $\til{M}$.
For each $\sigma \in \supp (C_k)$, choose a lift $\til{\sigma}$ of $\sigma$
whose center of gravity is in $D$. Since the bundle of frames over $D$
is compact, there is some ideal tetrahedron $\Delta'$ in $\H^n$ with center
of mass in $D$ for which some definite mass of $\til{\sigma}$ is geometrically
close to $\Delta'$ on a big set. Call $S$ the set of lifts sufficiently
close to $\Delta'$. It follows
that $\partial S$ is geometrically close
to $\partial \Delta$. Since $\partial \til{C_k} = 0$, most of the mass of
this boundary must be absorbed by simplices which are geometrically close to
being regular and ideal.

It follows that for $\delta$ sufficiently small, we can find lifts
$\til{\sigma}$ whose images are close on a big set
to the ideal simplices obtained by reflecting
$\Delta'$ in each of its boundary faces. Let $\alpha_k(\Delta) = \Delta'$.
Continuing inductively, if we propagate outwards in $X$ until we 
cover a big ball, we can find corresponding simplices 
in $\supp(\til{C_k})$ which agree with the simplices in $\alpha_k\text{dev}(T)$
to within a suitable tolerance, by taking $\delta$ sufficiently small.
\end{pf}

\begin{cor}[Jungreis]\label{simplices_everywhere}
For a hyperbolic $M^n$ with $n\ge 3$, a regular ideal simplex $\Delta \subset \H^n$
and any sequence $C_j$ of geodesic chains whose norms converge to $\| M\|_G$,  
for sufficiently large $j$ there is a $\sigma_j$ in the support of $C_j$ which
lifts to a geodesic simplex arbitrarily close to $\Delta$.
\end{cor}
\begin{pf}
By the previous lemma, there are simplices in the support of $\til{C_j}$
which stay close to $\alpha_k\text{dev}(T)$ on a big ball about
some fixed point. We can identify the set of framed ideal regular simplices in
$\H^n$ with $Isom^+(\H^n)$. In the limit, the set of ideal regular simplices in
the support of $\til{C_j}$ is invariant under the action of 
$\pi_1(M)$ on the left and $\text{dev}(Aut(T))$ on the right.
If $n>3$, $\text{dev}(Aut(T))$ is already dense in $Isom^+(\H^n)$. In dimension $3$,
following Jungreis and Ratner, using both the left and right actions we can find
simplices in the support of $C_j$ arbitrarily close to $\Delta$ for $j$
sufficiently large. See \cite{dJ97} for details.
\end{pf}

\begin{defn}
Say that a foliation $\F$ of a hyperbolic $n$--manifold is 
{\em asymptotically separated} if for some leaf
$\lambda$ of $\til{\F}$, there are a pair of open hemispheres
$H^+,H^- \subset \H^n$ in the complement of $\lambda$ which are
separated by $\lambda$.
\end{defn}

\begin{exa}
For $\F$ a finite depth foliation which is not a perturbation of a
surface bundle over a circle, the compact leaves lift to quasigeodesically
embedded planes in $\H^3$. Hence every leaf has the separation property,
and $\F$ is asymptotically separated.
\end{exa}

\begin{thm}\label{asymptotically_separated}
Suppose $M^n$ is hyperbolic and $\F$ is asymptotically separated. 
Then $$\| [M]\|_G < \| ([M],\F)\|_{FG}$$
\end{thm}
\begin{pf} By passing to a finite cover if necessary, we can assume that
$\F$ is co--oriented.
If we can show that every chain whose norm is sufficiently 
close to the Gromov norm contains an edge whose endpoints are not joined 
by an arc transverse to $\F$, then we will be done. 

Let $\lambda$ be a leaf of $\til{\F}$ 
and $H^+,H^-$ a pair of hemispheres in the complement
of $\lambda$ as provided by the definition of asymptotically separated. These
determine a pair of disks $D^+,D^- \subset S^{n-1}_\infty$ above and below
$\lambda$ respectively. Let $\alpha \in \pi_1(M)$ be an element taking the
complement of $D^+$ inside $D^+$. Then any infinite line from 
$D^-$ to $\alpha(D^-)$ must fail to be transverse to $\til{\F}$ somewhere,
since when it crosses $\lambda$ it is going in the positive
direction, and when it crosses $\alpha(\lambda)$ it is going in the
negative direction, with respect to the co--orientation on $\til{\F}$
which is preserved by $\alpha$.

It is easy to find an ideal regular simplex
$\Delta$ which has a pair of endpoints in $D^-$ and $\alpha(D^-)$ 
respectively. For any chain $C_j$ with norm sufficiently close to $\| [M]\|_G$,
there is a $\sigma$ in the support of $C_j$ whose geodesic representative
stays very close to $\Delta$ on an arbitrarily large compact piece. Such a
$\sigma$ has endpoints on incomparable leaves, and therefore cannot be
straightened (keeping its endpoints fixed) to a transverse simplex. It follows
that no transverse chain can have norm too close to $\| [M]\|_G$, and the
strict inequality is proved.
\end{pf}

\subsection{Limit sets of leaves of taut foliations}

\vskip 12pt

To investigate the asymptotic separation property, we must investigate the
limit sets of leaves of taut foliations.

\begin{lem}\label{above_below}
Let $\F$ be taut, and let $\lambda$ be a leaf of $\til{\F}$ on which we have
chosen a co--orientation. Denote by
$\lambda_\infty$ the limit set of $\lambda$. Then each region $D$ in the
complement of $\lambda_\infty$ is either above or below $\lambda$, in the
sense that for any two sequences $\lbrace p_i \rbrace \subset \H^3$ and
$\lbrace q_i \rbrace \subset \H^3$ with $p_i \to p \in D$ and 
$q_i \to q \in D$, the points $p_i,q_i$ are eventually on the same side of
$\lambda$.
\end{lem}
\begin{pf}
The points $p$ and $q$ can be joined by an arc $\alpha$ in $S^2_\infty$ which
avoids $\lambda_\infty$. This arc $\alpha$ is the Hausdorff limit in
$\bar{\H^3}$ of a sequence of arcs $\alpha_i$ joining $p_i$ to $q_i$. If each
of the $\alpha_i$ intersected $\lambda$, this would give rise to a sequence
of points in $\lambda$ converging to some point in $\alpha$, contrary to the
hypothesis that $\alpha$ avoids $\lambda_\infty$. It follows that
$p_i$ and $q_i$ are eventually on the same side of $\lambda$, and therefore
the ``side'' of $D$ is unambiguously defined.
\end{pf}

Notice that ``above'' and ``below'' as defined in the previous lemma are
{\em not} the same as $<$ and $>$ in the partial order on $L$. Each leaf
in the universal cover of a taut foliation has two sides; a co--orientation
on the leaf defines one of the sides to be above and one below, and {\em every}
other leaf falls into one of these two possibilities. This does {\em not} 
define a partial ordering on leaves.

\begin{defn}
Say that a foliation $\F$ has {\em two--sided branching} if in the
partial order on the leaf space $L$ of $\til{\F}$, there are triples of
leaves $\lambda,\lambda^+_l,\lambda^+_r$ and $\mu,\mu^-_l,\mu^-_r$ such that
$$\lambda < \lambda^+_l,\lambda < \lambda^+_r$$
$$\lambda^+_l \text{ and } \lambda^+_r \text{ are incomparable}$$
$$\mu^-_l < \mu, \mu^-_r < \mu$$
$$\mu^-_l \text{ and } \mu^-_r \text{ are incomparable}$$
\end{defn}

Observe that if $\F$ is taut and has two--sided branching, then we may
choose any leaf as $\mu = \lambda$. Moreover, if $\F$ is not co--orientable,
or covers some foliation which is not co--orientable, then either $\F$ is
$\R$--covered or it has two--sided branching.

\begin{exa}
Let $\F$ be a foliation of $T^3$ with one horizontal torus leaf, and the
complementary $T^2 \times I$ foliated as a Reeb foliation of the annulus
$\times S^1$. Then a pair of transversals whose initial segments agree
and cross the horizontal torus leaf must thereafter be leafwise homotopic;
that is, there is no branching in the positive direction from that point on.
Thus, we cannot choose $\mu = \lambda$ in this example.
\end{exa}

\begin{exa}\label{oneway_branching}
In \cite{gM91} G. Meigniez constructs examples of taut foliations which
branch on only one side, say the negative side. 
Furthermore, some of these examples are
obtained as perturbations of surface bundles over circles, and therefore
have pseudo--Anosov flows transverse to them. In \cite{dC99b} we construct
new examples of such foliations, and show that this situation holds in
general: taut foliations of atoroidal $3$--manifolds with one--sided branching
have transverse pseudo--Anosov flows which are {\em regulating}: that is,
flow lines are properly embedded in the leaf space of the universal cover.

Easy examples of foliations with one--sided branching are obtained by
starting with $\R$--covered foliations with (approximately) 
projectively invariant transverse measures, and then taking branched covers
over a curve which lifts to a line in $\til{M}$, one end of which is properly
embedded in the leaf space and one end of which is not.
\end{exa}

\begin{thm}
Let $\F$ be a taut foliation. If $\F$ has two--sided branching and
is not asymptotically separated, then
every leaf $\lambda$ of $\til{\F}$ has $m(\lambda_\infty) > 0$, where
$m$ is some normalized Lebesgue measure on the sphere at infinity 
$S^2_\infty$ of $\H^3$. In fact, $\lambda_\infty$ must have non--empty
interior.
\end{thm}
\begin{pf}
Assume without loss of generality that $\F$ is co--oriented.

Suppose that $m(\lambda_\infty) = 0$ for some $\lambda$. Then certainly
there is some complementary domain to $\lambda_\infty$ in $S^2_\infty$.
If there are domains $D^\pm$ both above and below $\lambda$, in the sense of
lemma~\ref{above_below} then there are half--spaces $H^\pm$ bounded by
circles in $D^\pm$ which avoid $\lambda$, and $\F$ is asymptotically 
separated. Otherwise without loss of generality, all the complementary
regions to $\lambda_\infty$ are contained above $\lambda$. It follows
that the subset of $\til{M}$ below 
$\lambda$ has limit set exactly equal to $\lambda_\infty$.

If $\F$ is $\R$--covered, one knows that $\lambda_\infty = S^2_\infty$
for every $\lambda$, so we may assume $\F$ is not $\R$--covered. (see
e.g. \cite{sF92})

If $\F$ has two--sided branching, then
there are a pair of positive transversals to $\til{\F}$ emanating from
$\lambda$ and ending on two incomparable translates $\alpha(\lambda),
\beta(\lambda)$ both $> \lambda$ in the partial order on $L$.
Now, the subset of $\til{M}$ below $\alpha(\lambda)$ has limit set equal to
$\alpha(\lambda_\infty)$, which has measure $0$. However, the subset of
$\til{M}$ above $\beta(\lambda)$ is itself a subset of the subset of $\til{M}$
below $\alpha(\lambda)$. It follows that we can write $\H^3$ as the union of
two sets (the sets above and below $\beta(\lambda)$), each of which has
a limit set of measure $0$, which is absurd. More generally, if
$\lambda_\infty$ has no interior, we could write $S^2_\infty$ as the
union of two closed sets without interior, which is absurd.
\end{pf}

The following theorem is proved in \cite{sF98}:

\begin{thm}[Fenley]
Let $\F$ be a Reebless foliation in $M^3$ closed, hyperbolic. Suppose
that $\lambda_\infty \ne S^2_\infty$ for some $\lambda$, and assume that
there is branching in the positive and negative directions of $\til{\F}$.
Then there is a $k<2$ such that the limit set of every leaf has Hausdorff
dimension less than $k$. In particular, every such limit set has zero
Lebesgue measure.
\end{thm}

\begin{cor}
No leaf in the universal cover of a taut foliation of a hyperbolic
$3$--manifold with two--sided branching 
is quasi--isometric (as a subset of $\H^3$) to a totally
degenerate surface group.
\end{cor}
\begin{pf}
The limit set of a totally degenerate surface group is a dendrite: a closed
set of measure $0$ whose complement is connected 
(see e.g. \cite{cMcM96}). If the limit set of
a leaf of a taut foliation has measure $0$, it has at least $2$ complementary
regions: one above and one below.
\end{pf}

\begin{cor}
Let $\F$ have two--sided branching. Then either $\lambda_\infty = S^2_\infty$
for every leaf $\lambda$ of $\til{\F}$ or every leaf is asymptotically
separated, and the foliated norm of $[M]$ is strictly greater than the usual
norm.
\end{cor}

Fenley has conjectured that $\lambda_\infty = S^2_\infty$ iff $\F$ is 
$\R$--covered.

By contrast, we have the following theorem:

\begin{thm}
Suppose that $\F$ is a taut foliation with one--sided branching. 
Then there is an equality of norms $$\|[M]\|_G = \| ([M],\F) \|_{FG}$$
\end{thm}
\begin{pf}
We declare that the branching takes place in the positive direction, with
respect to some choice of co--orientation on $\til{\F}$.

Lift the singular maps in the support of a chain $C = \sum_i r_i \sigma_i$ to
maps $\til{\sigma_i}:\Delta^3 \to \til{M}$. It is possible that some
vertices of $\Delta^3$ are mapped by some $\til{\sigma_i}$ to incomparable
leaves of $\til{\F}$. However, any pair of points in $\til{M}$ can be made
comparable after a finite isotopy in the negative direction. Moreover,
if two points are already on comparable leaves, then they are still on
comparable leaves after such an isotopy. 
Since there are only finitely many $\til{\sigma_i}$, we
can push the images of the vertices under ${\sigma_i}$ 
in the negative direction to get a new chain,
homotopic to $C$, for which each $\til{\sigma_i}$ sends the vertices
of $\Delta^3$ to comparable leaves of $\til{\F}$. By 
lemma~\ref{endpoints_ordered}, we can straighten this new chain relative to its
vertices to be transverse to $\F$.
\end{pf}

\begin{cor}
If $\F$ has one--sided branching, the leaves of $\til{\F}$ are not
asymptotically separated.
\end{cor}

\subsection{Foliations with Reeb components}

\vskip 12pt

One might suppose at least that the existence of Reeb components should be
detected by the foliated Gromov norm. However, this is not the case, as the
following example shows.

\begin{exa}
Let $M = S^2 \times S^1$ and $\F$ the standard foliation by two Reeb 
components. Let $C$ be any
chain representing $M$. Let $\pi:M \to M$ be the unique connected double 
cover of $M$. Then $\pi^*\F = \F$ up to isotopy. However, $\pi_* C = 2[M]$,
so the sequence $2^{-n} \pi_*^n C$ of chains can be made transverse after
isotopy and have norm $\to 0$. Hence $$\| [M] \|_G = \| ([M],\F)\|_{FG} = 0$$ 
in this case.
\end{exa}

Despite this example, it is easy to arrange a sequence of foliations
of a given manifold $M$ with more and more Reeb components where the
foliated Gromov norm grows without bound, as the following theorem shows.

\begin{thm}
There is a function $k:\N \to \R$ with
$\lim_{n \to \infty} k(n) = \infty$ such that if $\F$ is any foliation of
a $3$--manifold $M$ with $n$ generalized Reeb components whose complement 
is atoroidal, then $$\| ([M],\F)\|_{FG} \ge k(n)$$
\end{thm}
\begin{pf}
A generalized Reeb component, also known as a dead end component, is
a region of the foliated manifold bounded by torus or Klein bottle compact
leaves, such that no path transverse to the foliation which enters
the component can leave again.

Let $C = \sum_i r_i \sigma_i$ be a transverse chain representing $[M]$ and
let $N \subset M$ be a dead end component. Suppose that $\sigma_i$ is a
singular simplex whose image intersects $N$. Let $p \in \Delta^3$ be such
that $\sigma_i(p) \in N$ and let $\alpha$ be a path in $\Delta^3$ transverse
to the foliation induced by $\sigma_i^{-1}(\F)$ running from the top to the
bottom vertices which passes through $p$. Then $\sigma_i(\alpha)$ is 
transverse to $\F$ and intersects $N$; it follows that the image of at least
one of the vertices of $\Delta^3$ must be contained in $N$.

If we ``truncate'' $M$ by removing the Reeb components, we get a $3$--manifold
$M'$ with at least $n$ torus or Klein bottle cusps. By the fact above, each
truncated simplex can be collapsed to an edge or a face, or else 
is a normal simplex possibly with some ideal points. 
The resulting truncated chain $C$ represents 
$[M'] \in H_3(M',\partial M';\R)$, so the
foliated Gromov norm of $M$ can be estimated by the usual Gromov norm of
$M'$. Then set $k(n)$ equal to the minimum Gromov norm of a hyperbolic
$3$--manifold with $n$ cusps.

After the work of Thurston (see e.g. \cite{rBcP92}),
one knows that $\lim_{n \to \infty} k(n) \to \infty$.
\end{pf}

\section{Extending the norm to $H_*(M;\R)$}
\subsection{Semicontinuity of the norm}

\vskip 12pt

It is clear that the definition of the foliated Gromov norm can be extended
to a norm on $H_i(M;\R)$ for a manifold $M$ foliated by $\F$. As before,
for each homology class $\mu$ we can consider transverse
singular chains representing $\mu$, and take the $L_1$ norm of such
representatives. Denote the value of this norm on a class $\mu$ by
$\| (\mu,\F)\|_{FG}$. 

Notice that unlike the usual Gromov norm, this norm may be non--trivial 
even on $H_1(M;\R)$, as the following example shows:

\begin{exa}
Let $\F$ be the foliation of $T^2 \times I$ obtained by multiplying a
Reeb foliation of the cylinder $S^1 \times I$ by $S^1$. Glue the top and
bottom of $T^2 \times I$ together to get a foliation of $T^3$ also denoted
by $\F$. Let $\alpha \in H_1(M;\Z)$ be the generator obtained from
the $I$ factor by the gluing. Let $\beta \in H_1(M;\R) = r\alpha$. Then
$\|(\beta,\F)\|_{FG} \ge r/2$. For, each $\til{\sigma}:\Delta^1 \to \R^3$
obtained by lifting a map in the support of a chain representing
$\beta$ must have length $\le 2$ in its projection to the vertical
factor, since such a chain cannot cross the torus leaf twice.
\end{exa}

\begin{thm}\label{semicontinuous}
Let $\F_i$ be a sequence of taut foliations of $M^n$ which converge geometrically
(as $(n-1)$--plane fields) to $\F$. Then
$$\| (A,\F)\|_{FG} \ge \limsup \| (A,\F_i)\|_{FG}$$
for any $A \in H_*(M;\R)$.
\end{thm}
\begin{pf}
Let $C_j = \sum_i r_{ij} \sigma_{ij}$ be a sequence of cycles transverse to
$\F$ representing $A$ whose norms converge to $\| (A,\F)\|_{FG}$. Then
for any $j$, every $\sigma_{ij}$ is transverse to $\F$ and therefore by
compactness there is an $\epsilon_j$ such that the $1$--skeleta of the
images of $\Delta^n$ under the $\sigma_{ij}$ make an angle of at least
$\epsilon$ with $\F$ everywhere. It follows that for sufficiently large
$k$, the $1$--skeleta of $\supp(C_j)$ is transverse to $\F_k$. By
lemma~\ref{endpoints_ordered} we can straighten $C_j$ to be transverse
to $\F_k$.
\end{pf}

This implies that the foliated Gromov norm is lower semi--continuous: the
norm can jump up at a limit, but never down. The following example shows,
however, that the norm is not actually continuous.

\begin{exa}
Suppose $M^3$ is hyperbolic and fibers over the circle, and has $b_2 \ge 2$. Let
$\F_i$ be a sequence of fiberings contained in some top dimensional face of
the Thurston norm converging to some foliation $\F$ which is at a vertex.
Then $\|([M],\F_i)\|_{FG} = \| [M]\|_G$ by theorem~\ref{R_covered}. On the
other hand, $\F$ contains a quasigeodesically embedded compact leaf, so
$\|([M],\F)\|_{FG} > \| [M]\|_G$ by theorem~\ref{asymptotically_separated}.
\end{exa}

An interesting phenomenon in the theory of foliations occurs when a sequence
of isotopies of a fixed foliation $\F$ converges geometrically to a 
topologically distinct foliation $\G$. We give a simple example of this
phenomenon.

\begin{exa}
Let $S$ be the cylinder  $I \times S^1$ foliated by horizontal circles
$\text{point} \times S^1$. For an end--preserving homeomorphism $f:I \to I$
we can produce a foliation $\F_f$ of $T^2 \times I$ which is the suspension
of the foliation of $S$ by the map $f \times \text{id}:S \to S$.
Any two topologically conjugate maps $f,g:I \to I$ give isotopic foliations.
Now, it is well--known that any two strictly increasing homeomorphisms of the
open interval to itself are topologically conjugate. One can easily find a
sequence $f_i$ of these which converge (as maps $I \to I$) to the identity.
The foliations $\F_{f_i}$ are all isotopic, but distinct from $\F_{\text{id}}$.
\end{exa}

In any case, upper--semicontinuity of the norm implies the following 

\begin{cor}
Let $\F,\G$ be taut foliations of $M$ and suppose 
$$\|(A,\F)\|_{FG} > \|(A,\G)\|_{FG}$$ for some $A$. 
Then no sequence of isotopies of $\F$ can converge geometrically to $\G$.
\end{cor}

\vskip 12pt

Dual to the $L_1$ norm on $C_*(M;\R)$ defined by a foliation, there is
an $L_\infty$ norm on $C^*(M;\R)$ defined as the supremum of the value of
the cochain on transverse singular maps. There is an associated
{\em foliated bounded cohomology}, denoted by $H^*_{\F}(M;\R)$.
This may contain nontrivial elements even in dimension $1$, in contrast
with the usual bounded cohomology.

\subsection{Length of a free homotopy class}

There is a homotopy--theoretic refinement of the norm on $H_1$.
Say that a free homotopy class $[\alpha]$ of loops is {\em transverse}
to $\F$ if $\alpha$ is freely homotopic to a
transverse circle. More generally, define the {\em length} of
$[\alpha]$, denoted $\ell([\alpha])$, 
to be the minimum number of subdivisions of $S^1$ needed to
make a representative transverse on each subdivision. Say that
this length is $0$ if no subdivision is necessary: that is, if
some loop representing $\alpha$ is either transverse as a circle to
$\F$ or can be homotoped into a leaf of $\F$.

Notice that for a co--oriented foliation, $\ell$ takes on only even values.

\begin{lem}
The length of a free homotopy class of loops $[\alpha]$
 is upper semi--continuous in the
geometric topology.
\end{lem}
\begin{pf}
The only non--obvious point to check is that if $[\alpha]$ has a representative
which is contained in a leaf of $\F$, then $\ell([\alpha])=0$ for any
$\G_i$ sufficiently close to $\F$ in the geometric topology.
Lift $\alpha$ to an arc $\til{\alpha}$ in $\til{M}$ contained in a leaf
of $\til{\F}$. Then there is a big ball $B$ containing $\til{\alpha}$ which
is foliated in a standard way by $\til{\F}$. Therefore, for any $\G_i$
sufficiently close to $\F$ in the geometric topology, $\til{\G_i}$ foliates
some slightly smaller ball, also containing $\alpha$, in a standard way.
The arc $\alpha$ may be made transverse or tangential to $\til{\G_i}$ in
this ball, implying that $\ell([\alpha])=0$ for $\G_i$.
\end{pf}

\begin{lem}
For any co--oriented taut foliation $\F$ which has two--sided branching, 
there is an $[\alpha] \in \pi_1(M)$ with $\ell([\alpha])\ge 2$.
\end{lem}
\begin{pf}
Let $\tau_1,\tau_2$ be two positive transversals to $\til{\F}$
emanating from the same point whose upper endpoints are on incomparable
leaves. Let $\sigma_1,\sigma_2$ be two negative transversals to $\til{\F}$
emanating from the same point whose lower endpoints are on incomparable
leaves.

Map $\tau_1 \cup \tau_2$ and $\sigma_1 \cup \sigma_2$ 
downstairs to $M$. Since $\F$ is taut, one can extend the image of
$\tau_1$ in the positive direction by an arc $\rho_1$
until it joins up with $\sigma_1$, and
do the same with $\tau_2$. The union makes up a loop
$\alpha$, consisting of two transverse arcs 
$\alpha_1 = \tau_1 \cup \rho_1 \cup \sigma_1$
and $\alpha_2 = \tau_2 \cup \rho_2 \cup \sigma_2$.
We orient $\alpha$ so that $\alpha_1$ is positive and $\alpha_2$
is negative.
Let $\til{\alpha}$ be a lift to $\til{M}$, and consider its projection
to the leaf space $L$ of $\til{\F}$: this consists of an alternating
sequence of positive and negative arcs, which pass over a branch of $L$
at each stage. If $t_i$ and $b_i$ denote the alternating sequence
of top and bottom leaves of the projection, then for all $i \in \Z$,
$$t_i > b_i, t_i > b_{i-1}$$
$$b_i \text{ and } b_{i+1} \text{ are incomparable}$$
$$t_i \text{ and } t_{i+1} \text{ are incomparable}$$
We can find a sequence of points $m_i$ with
$b_i < m_i < t_i$ and each $m_i,b_{i-1}$ and $m_i,t_{i+1}$ pairwise incomparable.
It is clear that for any
$\alpha'$ homotopic to $\alpha$, the projection of the corresponding
lift $\til{\alpha'}$ to $L$ must intersect each $m_i$. In particular,
such a lift intersects incomparable leaves of $\til{\F}$, so $\alpha'$ cannot
be transverse. See figure 3.
\end{pf}

\begin{figure}[h]
\scalebox{.5}{\includegraphics{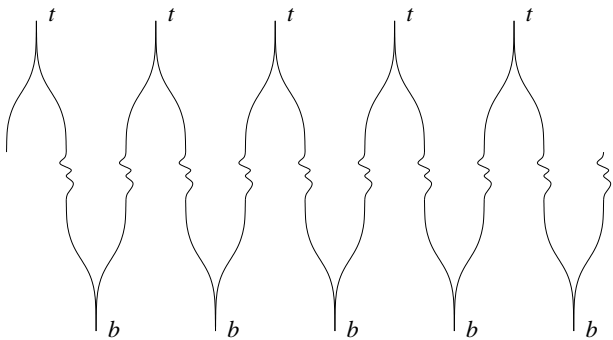}}
\caption{The sequence of lifts of $\alpha_1$ and $\alpha_2$, projected to
$L$ alternately branches in the positive and
negative direction; the same is true of any homotopic lift.}
\end{figure}

The following corollary answers a question posed by W. Thurston:

\begin{cor}
Let $\F$ with branching on at most one side and 
$\G$ with two--sided branching, be taut foliations of
$M$. Then there is no sequence of isotopies $\G_i$ of $\G$ which converges
geometrically to $\F$.
\end{cor}
\begin{pf}
Lift to a finite cover where the foliations are co--oriented. Then
observe that for a foliation with branching on at most one side, 
the length of any free homotopy class is $0$. For, suppose $\F$ does
not branch in the positive direction, and let $\alpha$ be any loop in
$M$, and suppose $\alpha$ is in general position with respect to
$\F$ so that it has a finite number of isolated minima and maxima. Let
$p$ be such a minimum, lying between maxima $q,r$. Lift the segment
$\tau$ between $q,r$ to $\til{\tau}$ in $\til{M}$. Then $\til{q},\til{r}$
are comparable, since they are both $> \til{p}$ and by hypothesis there is
no branching in the positive direction, so without loss of generality
we can assume $\til{q} \le \til{r}$. So we can push the local minima
corresponding to $p$ in the positive direction until it cancels the local maxima
corresponding to $q$ without introducing any new critical points; in particular,
the number of critical points can be reduced. Continuing inductively, they can all
be eliminated.
\end{pf}

\subsection{Virtually fine triangulations}

\begin{defn}
Let $M$ be a $3$--manifold. A triangulation $\tau$ of $M$ is
{\em fine} if for every taut foliation $\F$ of $M$ $\tau$ can be
isotoped to be transverse to $\F$. A triangulation $\tau$ is
{\em virtually fine} if for every taut foliation $\F$ of $M$ there is
a finite cover $\widehat{M}$ of $M$ such that the pulled--back triangulation
$\widehat{\tau}$ can be isotoped to be transverse to the pulled--back
foliation $\widehat{\F}$.
\end{defn}

One of the main theorems of \cite{dG99} states that for any $M$ there is
a fine triangulation $\tau$. This leads naturally to the question of
what conditions are necessary and sufficient on a triangulation to be
fine. An obvious condition is that the triangulation admit a transverse
foliation {\em locally}. For a geodesic triangulation of a hyperbolic
manifold, this is obvious, since in the projective model of hyperbolic
space, a hyperbolically geodesic triangulation looks like a 
Euclidean triangulation of the ball,
and a foliation by horizontal planes will be transverse.
It has been an open question whether {\em every}
geodesic triangulation of a hyperbolic $3$--manifold $M$ is
virtually fine.

It turns out this guess is incorrect: there are geodesic triangulations whose
simplices have diameters arbitrarily small compared to the injectivity radius
of the ambient manifold, which cannot be made transverse to certain
taut foliations.

\begin{thm}
Let $M$ be a hyperbolic $3$--manifold, and $\F$ {\em any} taut foliation with
$2$--sided branching. Then there is a geodesic triangulation $\tau$ of
$M$ which cannot be made transverse to $\F$.
Furthermore, $\tau$ cannot be made transverse to $\F$ in any finite cover
(i.e. $\tau$ is not virtually fine).
\end{thm}
\begin{pf}
Let $\gamma$ be a closed loop with $\ell(d) = l \ge 2$. Then we can choose a
geodesic representative of $\gamma$ and make it an edge of a geodesic
triangulation. Such a triangulation can obviously not be made transverse
to $\F$. Now, for any finite cover $\widehat{M}$ of $M$, we can lift
$\gamma$ to some $\widehat{\gamma}$ which covers $\gamma$ with degree $d$, and
has $d$ segments in the lifted triangulation.
With respect to $\widehat{\F}$, the new length of $\widehat{\gamma}$ is
$ld$, so the lifted triangulation cannot be made transverse to the lifted
foliation.
\end{pf}

\section{Laminations and order trees}

\subsection{Genuine laminations}

Laminations of $3$--manifolds are defined in \cite{dGuO89}.

\begin{defn}
A {\em lamination} in a $3$--manifold is a foliation of a closed subset of $M$
by $2$ dimensional leaves. The complement of this closed subset falls into
connected components, called {\em complementary regions}.
A lamination is {\em essential} if it contains
no spherical leaf or torus leaf bounding a solid torus, and furthermore if
$C$ is the closure (with respect to the path metric in $M$)
of a complementary region, then $C$ is irreducible and
$\partial C$ is both incompressible and {\em end incompressible} in $C$.
Here an end compressing disk is a properly embedded $(D^2 - (\text{closed
arc in } \partial D^2))$ in $C$ which is not properly isotopic relative
to the $\partial$
in $C$ to an embedding into a leaf. Finally, an essential lamination is
{\em genuine} if it has some complementary region which is not an $I$-bundle.
\end{defn}

An essential lamination simultaneously generalizes both Reebless foliations
and incompressible surfaces. It is {\em not} true that an essential
lamination lifts in a finite cover to a co--orientable lamination.
Consequently the leaf space of an essential lamination in the universal
cover does not carry a natural partial order. The leaf space of a foliation
is like a train--track: there is a natural combing near any branch point.
The leaf space of a lamination is more like a tree: there is no natural 
way to say whether branches approach a branch point from the same or from
opposite directions.

Nevertheless, we can still
talk about transversality of a simplicial map in a laminated context.

\begin{defn}
Let $\Lambda$ be an essential lamination of $M$.
A map $\sigma:\Delta^1 \to M$ is {\em transverse} if there is no
{\em back--tracking}; i.e. there is no subinterval of $\Delta^1$ whose image
can be homotoped relative to its endpoints into a leaf of $\Lambda$.
A map $\sigma:\Delta^i \to M$ with $i\ge 2$ 
is {\em transverse} if the induced lamination
$\sigma^{-1}(\Lambda)$ of $\Delta^i$ is non--singular and can be completed
to an affine foliation of $\Delta^i$
\end{defn}

Denote by $\|(A,\Lambda)\|_{LG}$ the norm of a homology class $A$ with
respect to $\Lambda$.

For $\sigma:\Delta^3 \to M$ a transverse map
and $\Lambda$ nowhere dense, we can perturb
$\sigma$ to be nondegenerate; that is, so that $\sigma^{-1}(\Lambda)$ is a collection of
normal triangles and quadrilaterals compatible with a total ordering
of the vertices of $\Delta^3$. This can be done by wiggling $\sigma$ slightly
so that no vertex is taken into a leaf of $\Lambda$.

Let $\T$ be a triangulation of $M$, and let
$n(\T)$ denote the number of tetrahedra in $\T$. 
Then any {\em minimal} genuine lamination
$\Lambda$ (i.e. one with every leaf dense in $\Lambda$) can be put in normal
form with respect to $\T$, by a theorem of M. Brittenham
\cite{mB95}. For instance, an incompressible surface is an example of a
minimal lamination. This establishes the following estimate

\begin{thm}
Let $\Lambda$ be a minimal genuine lamination of $M$. Then
$$\|([M],\Lambda)\|_{LG} \le \min_{\T} 4n(\T)$$
\end{thm}
\begin{pf}
Let $\T$ be a triangulation of $M$. Then we can isotope $\Lambda$ to be
in normal form with respect to $\T$. Now we can subdivide $\T$, replacing
each tetrahedron by $4$ tetrahedra, each of which only contains normal
disks compatible with a total ordering on its vertices.
\end{pf}

On the other hand, corollary~\ref{simplices_everywhere} implies

\begin{thm}
For $S$ an incompressible surface in a hyperbolic $3$--manifold $M$, 
either $S$ is a fiber of a fibration of $M$ over $S^1$ or
$$\|([M],S)\|_{LG} > \|[M]\|_{G}$$
\end{thm}
\begin{pf}
Either $S$ is a fiber of a fibration over $S^1$, or $S$ is quasigeodesic.
In the second case, we can find three lifts $\til{S}_1,\til{S}_2,\til{S}_3$
of $S$ to $\til{M} = \H^3$ which bound disjoint half--spaces. There is
an ideal triangle with one vertex in each of these half--spaces, and there
is a regular ideal tetrahedron, one of whose faces is this triangle. It
follows that any chain representing $[M]$ whose norm is sufficiently
close to $\|[M]\|_{G}$ cannot be transverse to $S$.
\end{pf}

\subsection{Order trees}

The following definition is found in \cite{dGuO89}.
 
\begin{defn}
An {\em order tree} is a set $T$ together with a collection $\S$ of
linearly ordered subsets called {\em segments}, each with distinct
least and greatest elements called the {\em initial} and {\em final}
ends. If $\sigma$ is a segment, $-\sigma$ denotes the same subset with
the reverse order, and is called the {\em inverse} of $\sigma$. The
following conditions should be satisfied:
\begin{itemize}
\item{$\sigma \in \S \Longrightarrow -\sigma \in \S$}
\item{Any closed subinterval of a segment is a segment (if it has more
than one element).}
\item{Any two elements of $T$ can be joined by a finite sequence of
segments $\sigma_i$ with the final end of $\sigma_i$ equal to the
initial end of $\sigma_{i+1}$.}
\item{Given a cyclic word $\sigma_0\sigma_1 \dots \sigma_{k-1}$ (subscripts
mod $k$) with the final end of $\sigma_i$ equal to the initial end of
$\sigma_{i+1}$, there is a subdivision of the $\sigma_i$ yielding a cyclic
word $\rho_0\rho_1 \dots \rho_{n-1}$ which becomes the trivial word when
adjacent inverse segments are canceled.}
\item{If $\sigma_1$ and $\sigma_2$ are segments whose intersection is a
single element which is the final element of $\sigma_1$ and the initial 
element of $\sigma_2$ then $\sigma_1 \cup \sigma_2$ is a segment containing
$\sigma_1$ and $\sigma_2$.}
\end{itemize}
\end{defn}

An order tree is topologized by the order topology on segments. We assume
in the sequel that our order trees are {\em $\R$--order trees} --- that is,
$2$nd countable order trees whose segments are order isomorphic to
compact intervals of $\R$.

Let $\Gamma = \pi_1(M)$ act by automorphisms on an order tree $T$
and suppose that we have a
$\Gamma$--equivariant surjective map $$\phi:\til{M} \to T$$

A singular map $\sigma:\Delta^i \to M$ is {\em transverse} if for any
lift $\til{\sigma}:\Delta^i \to \til{M}$ the composition 
$\phi\til{\sigma}:\Delta^i \to T$ maps $\Delta^i$ to a totally ordered
segment of $T$.

Say that a sequence of such maps $\phi_i:\til{M} \to T_i$ {\em converges}
to $\phi:\til{M} \to T$ if every map $\sigma:\Delta \to M$ transverse with
respect to $\phi$ is eventually transverse with respect to $\phi_i$, for
sufficiently large $i$.

For a representation $\rho:\Gamma \to \text{Aut}(T)$, an equivariant map
$\phi:\til{M} \to T$ and a homology class $\mu \in H_i(M;\R)$ we can
define a norm $\| (\mu,\phi)\|_{FG}$ as before as the $L_1$ norm on
the singular chains representing $\mu$ which are transverse with respect to
$\phi$. Observe that this norm does not really depend on the map
$\phi$, since it is determined up to equivariant homotopy by $\rho$, and
thus this is really a norm on $H_i(\Gamma;\R)$ depending only on $\rho$.

\begin{thm}
Let $\Gamma = \pi_1(M)$ and let $\phi_i:\til{M} \to T_i$ be a sequence
of equivariant maps for actions $\rho_i:\Gamma \to \text{Aut}(T)$.
Suppose this sequence converges to $\phi:\til{M} \to T$ equivariant for
some action $\rho:\Gamma \to \text{Aut}(T)$.
 
Then given $\mu \in H_i(M;\R)$ we have the inequality
$$\| (\mu,\phi) \|_{FG} \ge \lim \sup \|(\mu,\phi_j) \|_{FG}$$ 
\end{thm}
\begin{pf}
Any geometric chain in $M$ representing $\mu$ which is
transverse for $\phi$ will be transverse for $\phi_i$ for sufficiently
large $i$.
\end{pf}

Examples of group actions on order trees arise in the study of
essential laminations, where the lamination in the universal cover is
dual to an order tree which can be taken to be an $\R$--order tree by
replacing isolated leaves by foliated $I$--bundles over those leaves.
We have already seen from the example of foliations that this norm is
not trivial and can vary for different representations of a fixed group.
\vfill
\pagebreak

\end{document}